\documentstyle[11pt]{article}

\title{THE GALAXIES OF NONSTANDARD ENLARGEMENTS OF INFINITE AND TRANSFINITE GRAPHS: II}
\author{A. H. Zemanian}
\date{}

\begin{document}
\newcommand{\N} {I \kern -4.5pt N}
\newcommand{\R} {I \kern -4.5pt R}
\newcommand{\Z} {Z \kern -7.5pt Z}
\maketitle
\baselineskip21pt

{\ Abstract --- This report is an improvement of a prior report (Report 813).  
It sharpens the principal theorems (Theorems 4.2 and 11.2 of Report 813) 
while simplifying their proofs.  There are also several minor 
changes involving 
clarifications and corrections of misprints.

The Abstract of the prior report remains the same as follows:
The galaxies of the nonstandard enlargements of connected,
conventionally infinite graphs as well as of connected transfinite graphs are 
defined, analyzed, and illustrated by some examples.  
It is then shown that any 
such enlargement either has exactly one galaxy, its principal one, or 
it has infinitely many galaxies.  In the latter case, the galaxies are 
partially ordered by their ``closeness'' to the principal galaxy.
If an enlargement has a galaxy different from its principal galaxy,
then it has a two-way infinite sequence of galaxies that are 
totally ordered according to that ``closeness'' property.
There may be many such totally ordered seqences.

} 

\section{Introduction}

In this work we extend the idea of galaxies in the hyperreal line
$^{*}\!\R$ to nonstandard enlargements of conventionally infinite 
graphs and also of transfinite graphs.  We stipulate henceforth that 
every graph considered herein is connected.
Since graphs have structures 
much different from that of the real line $\R$, the enlargements 
of graphs have properties not possessed by $^{*}\! \R$.
The graphical galaxies of those enlargements comprise one
aspect of that distinctive complexity.  We will show that 
that any such enlargement has either one galaxy or infinitely many of them.
Moreover, just as $^{*}\! \R$ contains images of the real numbers, 
called the standard hyperreals, as well as hyperreals that are nonstandard, 
so too may the enlargement $^{*}\!G$ of 
a graph $G$ contain ``hypernodes,'' some of which are images 
of nodes of $G$ and others of which are nonstandard hypernodes.
In addition, there are ``hyperbranches'' incident to pairs of hypernodes;
some of these hyperbranches are images of branches of $G$, but there may be
others that are not.

The galaxies graphically partition $^{*}\!G$ in the sense the every
hypernode belongs to exactly one galaxy, and so too does every hyperbranch.
There is a unique galaxy, which we refer to as the ``principal 
galaxy,'' that contains the standard  hypernodes and possibly 
nonstandard hypernodes as well.  In the event that there are
infinitely many galaxies, those galaxies are partially
ordered according to how ``close'' they are 
to the principal galaxy.  In fact,  
if there is a galaxy different from the principal galaxy,
then there is a two-way infinite sequence of galaxies
that are totally ordered according to their ``closeness'' to
the principal galaxy.  There may be many such totally ordered sequences, but
a galaxy in one such sequence may not be comparable to a galaxy in 
another sequence according to that ``closeness'' property.

We speak of ``conventionally infinite'' graphs to 
distinguish them from transfinite graphs of ranks 1 or higher
\cite[Chapter 2]{tgen}, \cite[Chapter 2]{gn}.  Sections 2 through 4 
herein are devoted to the enlargements of conventionally infinite graphs.
The results for such enlargements extend to enlargements of 
transfinite graphs, but in more complicated ways.  We show this in 
Sections 5 through 11, but only for transfinite graphs of rank 1.
Results for transfinite graphs of still higher ranks 
are obtained similarly but in still more
complicated ways and with additional complexity in the symbols.
For the sake of brevity, the latter results are not included herein,
but they may be found in \cite{gal2} as well as in the archive
www.arxiv.org in the category ``mathematics'' under ``Zemanian.''

Our notations and terminology follow the usual conventions of nonstandard 
analysis.  $\N= \{0,1,2,\ldots\}$ is the set of 
natural numbers, and $^{*}\! \N$ is the set of hypernaturals.
The standard hypernaturals are (i.e., can be identified 
with) the natural numbers.  Also, $\langle a_{n}\rangle$ or 
$\langle a_{n}\!: n\in\N\rangle$ or $\langle a_{0},a_{1},
a_{2},\ldots\rangle$ denotes a sequence whose elements can be 
members of any set, such as the set $X$ of nodes in a 
conventional graph $G=\{X,B\}$, where $B$ is the set of branches,
a branch being a two-element set of nodes.  On the other hand, 
$[a_{n}]$ denotes an equivalence class of sequences, where
two sequences
$\langle a_{n}\rangle$ and $\langle b_{n}\rangle$ 
are taken to be equivalent if $\{n\!:a_{n}=b_{n}\}\in {\cal F}$,
where ${\cal F}$ is any chosen and fixed 
free ultrafilter.\footnote{Also called 
a nonprincipal ultrafilter.}  ${\cal F}$ will be so fixed throughout this work.
The $a_{n}$ appearing in $[a_{n}]$ are understood
to be the elements of any one of the sequences in the equivalence class.
At times, we will use the more specific notation 
$[\langle a_{0},a_{1},a_{2},\ldots\rangle]$.  More generally, 
we adhere to the notations and terminology appearing in \cite{go}.

The ordinals are denoted in the usual way: $\omega$ is the 
first transfinite ordinal.  With $\tau\in\N$, 
the product $\omega\cdot \tau$ is the sum of $\tau$ terms, 
each being $\omega$. 

\section{The Nonstandard Enlargement of a Graph}

Throughout Sections 2 to 4, we assume that the conventionally
infinite graph $G$ is connected and has infinitely many nodes.
The definition of a nonstandard graph that we use herein is given in 
\cite[Section 8.1]{gn}, a special case of which is the ``enlargement''
of a graph $G$.

Let us define the {\em enlargement} $^{*}\!G$ of $G$ 
here as well in order to remove 
any need for referring to \cite{gn}.  $G=\{X,B\}$ is 
now taken to be a conventional connected graph having an infinite set $X$
of nodes and therefore an infinite set of branches as well, each 
branch being a two-element set of nodes. Thus, there are no 
parallel branches (i.e., multiple branches). ${\cal F}$ will denote 
a chosen and fixed free ultrafilter.  
${\bf x}=[x_{n}]$ denotes an equivalence class of sequences 
of nodes as stated in the Introduction.  Specifically, 
$[x_{n}]$ and $[x'_{n}]$ are in the same equivalence class $\bf x$
if $\{n\!: x_{n}=x'_{n}\}\in{\cal F}$.
${\bf x}$ will be called a {\em 
hypernode}.\footnote{Our terminology should not be confused with that of 
a hypergraph---an entirely different concept \cite{be}.}
Thus, the set of all sequences of nodes from $G$ is partitioned 
into hypernodes.  $^{*}\!X$ denotes the set of hypernodes.  
If all the elements of one of the representative sequences 
$\langle x_{n}\rangle$ for a hypernode ${\bf x}=[x_{n}]$ are the same node
(i.e., $x_{n}=x$ for all $n$), then ${\bf x}=[x]$ can be 
identified with $x$;  in this case, ${\bf x}$ is called a {\em standard 
hypernode}.  Otherwise, 
${\bf x}=[x_{n}]$ is called a {\em nonstandard hypernode}.

We turn now to the definition of a ``hyperbranch.''  Let ${\bf x}=[x_{n}]$ and 
${\bf y}=[y_{n}]$ be two hypernodes.  Also, let ${\bf b}=[\{x_{n},y_{n}\}]$,
where $\langle \{x_{n},y_{n}\}\rangle$ is a sequence of pairs of nodes from 
$G$ such that, for almost all $n$, $\{ x_{n},y_{n}\}$
is a branch in $G$; that is, $\{n\!: \{x_{n},y_{n}\}\in B\}\in {\cal F}$.
It can be shown \cite[page 155]{gn} that this definition is independent
of the representative sequences $\langle x_{n}\rangle$ and
$\langle y_{n}\rangle$ chosen of ${\bf x}$ and ${\bf y}$ respectively and 
that we truly have an equivalence relation for the set of all 
sequences of branches from $G$. We let ${\bf b}=[\{x_{n},y_{n}\}]$ denote such 
an equivalence class and will call it a {\em hyperbranch};
we write ${\bf b}=\{{\bf x},{\bf y}\}$. Also,  $^{*}\!B$ will denote the set of all
hyperbranches.  If ${\bf x}=[x_{n}]$ and ${\bf y}=[y_{n}]$ are 
standard hypernodes, then ${\bf b}=[\{x,y\}]$ is called a 
{\em standard hyperbranch}.  Otherwise, ${\bf b}$ is called a 
{\em nonstandard hyperbranch}.

Finally, the pair $^{*}\!G=\{^{*}\!X,\,^{*}\!B\}$ denotes the 
{\em enlargement} of $G$.  It is a special case of a nonstandard graph,
as defined in \cite[page 155]{gn}.\footnote{If $G$ were a finite graph,
then every hypernode (resp. hyperbranch) could 
be identified with a node (resp. branch)in $G$, and $^{*}\!G$
would be identified with $G$.}

\section{Distances and Galaxies in Enlarged Graphs}

The {\em length} $|P_{x,y}|$ of any path $P_{x,y}$ connecting 
two nodes $x$ and $y$ in a graph $G$ is the number of branches in
$P_{x,y}$.  The {\em distance} $d(x,y)$ between $x$ and $y$ 
is $d(x,y)=\min\{|P_{x,y}|\}$, where the minimum is taken over all
paths terminating at $x$ and $y$.  In the trivial case, $d(x,x)=0$. 
$d$ satisfies the triangle inequality, namely, for any three nodes 
$x$, $y$, and $z$ in $G$, $d(x,y)\leq d(x,z)+d(z,y)$. In fact, $d$ 
satisfies the other metric axioms, too, 
and the set $X$ of nodes in $G$ along with $d$
is a metric space.

The metric $d$ can be extended into an internal function $\bf d$
mapping the Cartesian product $^{*}\! X\,\times \,^{*}\! X$ into the set of 
hypernaturals $^{*}\! \N$ as follows:  For any ${\bf x}=[ x_{n}]$ and 
${\bf y}=[y_{n}]$ in $^{*}\!X$, ${\bf d}$ is defined by 
\[ {\bf d}({\bf x},{\bf y})\;=\;[d(x_{n},y_{n}]\,\in\,^{*}\!\N. \]
By the transfer principle, we have, for any three hypernodes
${\bf x}$, ${\bf y}$, and ${\bf z}$,
\begin{equation}
{\bf d}({\bf x},{\bf z})\;\leq\; {\bf d}({\bf x},{\bf y})\,+\,{\bf d}(({\bf y},{\bf z}).  \label{3.1}
\end{equation}
From the point of view of an ultrapower construction, this means that
\[ \{n\!: d(x_{n},z_{n})\;\leq\;d(x_{n},y_{n})\,
+\,d(y_{n},z_{n}\}\;\in\;{\cal F}. \]
The other metric axioms, such as ${\bf d}({\bf x},{\bf x})=0$, are obviously
satisfied by ${\bf d}$. 

We define the ``galaxies'' of $^{*}\!G$ as nonstandard subgraphs
of $^{*}\!G$ by first defining the ``nodal galaxies.''  Two hypernodes
${\bf x}=[x_{n}]$ and ${\bf y}=[y_{n}]$ are taken to be in the same 
{\em nodal galaxy} $\dot{\Gamma}$ of $^{*}\!G$
if ${\bf d}({\bf x},{\bf y})$ is no 
greater that a standard hypernatural $\bf k$, that is, if there exists a natural 
number $k\in\N$ such that $\{n\!: d(x_{n},y_{n})\,\leq\,k\}\;\in\; {\cal F}$.
In this case, we say that ${\bf x}$ and ${\bf y}$ are {\em limitedly distant},
and we write ${\bf d}({\bf x},{\bf y})\leq {\bf k}$.

Let $N_{{\bf x},{\bf y}}$ be the set of all standard hypernaturals that are no
less than ${\bf d}({\bf x},{\bf y})$.  $N_{{\bf x},{\bf y}}$ is a 
well-ordered set, and 
therefore it has a minimum ${\bf k}_{{\bf x},{\bf y}}$.
So, we can say that ${\bf x}$ and ${\bf y}$ are in the same nodal galaxy
$\dot{\Gamma}$ if ${\bf d}({\bf x},{\bf y})={\bf k}_{{\bf x},{\bf y}}$.

{\bf Lemma 3.1.}  {\em The nodal galaxies partition the set $^{*}\!X$
of all hypernodes in $^{*}\!G$.}

{\bf Proof.}  The property of two hypernodes being limitedly distant
is a binary relation on $^{*}\!X$ that is obviously reflexive and symmetric.
Its transitivity follows directly from (\ref{3.1}).  Alternatively, 
we can use an ultrapower argument.  Assume that ${\bf x}=[x_{n}]$ and
${\bf y}=[z_{n}]$ are in some nodal galaxy and that ${\bf y}$ and
${\bf z}=[z_{n}]$ are in some nodal galaxy; we want 
to show that those galaxies are the same.  There exist
two standard natural numbers $k_{1}$ and $k_{2}$ such that 
$N_{{\bf x},{\bf y}}=\{n\!:d(x_{n},y_{n})\leq k_{1}\}\in{\cal F}$
and $N_{{\bf y},{\bf z}}=\{n\!:d(y_{n},z_{n})\leq k_{2}\}\in{\cal F}$.
Since $d(x_{n},z_{n})\,\leq\,d(x_{n},y_{n})+d(y_{n},z_{n})$,
\[ \{n\!: d(x_{n},z_{n})\leq k_{1}+k_{2}\}\;\supseteq\;
N_{{\bf x},{\bf y}}\cap N_{{\bf y},{\bf z}}\;\in\;{\cal F}. \]
So, the left-hand side is a set in ${\cal F}$.  Thus, ${\bf x}$ and ${\bf z}$ 
are limitedly distant, too, and ${\bf x}$, ${\bf y}$, and ${\bf z}$ are
all in the same nodal galaxy. $\Box$

We define a {\em galaxy} $\Gamma$ of $^{*}\!G$ as a maximal
nonstandard subgraph of  
$^{*}\!G$ whose hypernodes are all in the same nodal galaxy
$\dot{\Gamma}$;  that is, the hyperbranches of $\Gamma$
corresponding to $\dot{\Gamma}$ are all those pairs $\{{\bf x},{\bf y}\}$
such that ${\bf x},{\bf y} \in \dot{\Gamma}$.  We will say that a 
hypernode ${\bf x}$ is {\em in} $\Gamma$ when ${\bf x}\in\dot{\Gamma}$
and that a hyperbranch $\{{\bf x},{\bf y}\}$ is {\em in} $\Gamma$ when 
${\bf x},{\bf y}\in\dot{\Gamma}$.  
It follows from Lemma 3.1 that the galaxies of 
$^{*}\!G$ partition $^{*}\!G$ in the sense of graphical partitioning
(i.e., each hyperbranch is in one and only one galaxy).

The {\em principal galaxy} $\Gamma_{0}$ of $^{*}\!G$ is that unique
galaxy, each of whose hypernodes is limitedly distant from some
standard hypernode 
(and therefore from all standard hypernodes).
All the nodes in $G$ will be (i.e., can be identified with)
standard hypernodes in $\Gamma_{0}$, but there may be nonstandard
hypernodes in $\Gamma_{0}$ as well.  The following examples   
illustrate this point.

{\bf Example 3.2.}  Consider the endless (i.e., two-way infinite)
path:
\[ P\;=\;\langle\ldots,x_{-1},b_{-1},x_{0},b_{0},x_{1},b_{1}\ldots\rangle \]
with nodes $x_{k}$ and branches $b_{k}$, $k\in\Z$, $\Z$ 
being the set of integers.  The enlargement $^{*}\!P$ of $P$ has
hypernodes, each being represented by $[x_{k_{n}}]$ where
$\langle k_{n}\rangle$ is some sequence of integers.
Each hyperbranch is represented by  
$[\{x_{k_{n}},x_{k_{n}+1}]$.  The nodal galaxies are infinitely many 
because they correspond bijectively with the galaxies of the 
enlargement $^{*}\!\Z$ of $\Z$.  Moreover, the principal galaxy
$\Gamma_{0}$ of $^{*}\!P$ has only standard hypernodes and in fact is
(i.e., can be identified with) $P$ itself.  Also, every galaxy
is graphically isomorphic to $\Gamma_{0}$ and therefore to 
every other galaxy.  $\Box$

{\bf Example 3.3.} Now, consider a one-ended path:
\[ T\;=\;\langle x_{0},b_{0},x_{1},b_{1},x_{2},b_{2},\ldots\rangle  \]
Each hypernode in the enlargement $^{*}\!T$ of $T$ 
is represented by $[x_{k_{n}}]$, where $\langle k_{n}\rangle$
is some sequence of natural numbers.  Thus, $^{*}\!T$ has a hypernode set
$^{*}\!X$ that can be identified with the set $^{*}\!\N$ of 
hypernaturals.  Hence, $^{*}\!T$ has an infinitely of galaxies, too.
The principal galaxy $\Gamma_{0}$ of $^{*}\!T$ is the one-ended
path $T$.  However, any hypernode ${\bf x}=[x_{k_{n}}]$
in a galaxy $\Gamma$ different from $\Gamma_{0}$ will be such that, 
for every $m\in\N$,
$\{n\!: k_{n}>m\}\in{\cal F}$.  Such a hypernode is adjacent both to  
$[x_{k_{n}+1}]$ and to $[x_{k_{n}-1}]$, where we are free to replace
$x_{k_{n}-1}$ by, say, $x_{0}$ whenever $k_{n}=0$.
(The set $\{n\!: k_{n}=0\}$ will not be a member of ${\cal F}$ when 
${\bf x}=[x_{k_{n}}]$ is in $\Gamma$.)  Thus, ${\bf x}=[x_{k_{n}}]\in\Gamma$
has both a predecessor and a successor, which implies
that $\Gamma$ is graphically isomorphic to an endless path.  In fact, 
all the galaxies other than $\Gamma_{0}$ are isomorphic to each 
other, being identifiable with an endless path.  $\Box$

{\bf Example 3.4.}  Consider next the grounded, one-way infinite
ladder $L$ of Figure 1.  Now, for every $k\in\N$, 
$d(x_{k},x_{g})=d(x_{k},x_{k+1})=1$, and, for every 
$k,l\in\N$ with $|k-l|>1$, $d(x_{k},x_{l})=2$.  In this case, for 
every two hypernodes ${\bf x}$ and ${\bf y}$,
${\bf d}({\bf x},{\bf y})\leq[2]=2$.  Thus, every two hypernodes 
are limitedly distant 
from each other, which means that $^{*}\!L$ has only one galaxy, 
its principal galaxy $\Gamma_{0}$.  Now, $\Gamma_{0}$ has both 
standard and nonstandard hypernodes.  $\Box$

{\bf Example 3.5.}  Furthermore, consider the graph $G$ obtained 
from $L$ by appending a one-ended path $P$ starting at $x_{g}$, but
otherwise isolated from $L$, as shown in Figure 2.  
In this case, we again have an infinity of galaxies
by virtue of the isolation of $P$ from $L$.
The principal galaxy $\Gamma_{0}$ has both standard and nonstandard
hypernodes, its nonstandard hypernodes being due to $L$.
All the other galaxies are graphically isomorphic to  
an endless path (as in Example 3.3) and thus to each other, 
but not to $G$ and not to $\Gamma_{0}$. $\Box$

A subgraph $G_{s}$ of $G$ with the property that there exists a 
natural number $k$ such that $d(x,y)\leq k$ for all pairs of 
nodes $x,y$ in $G_{s}$ will be called a {\em finitely dispersed}
subgraph of $G$.  Example 3.5 suggests that the structures of the galaxies
other than $\Gamma_{0}$ do not depend upon any finitely dispersed 
subgraph of $G$.  This is true in general because the nodes 
$x_{n}$ in any representative $\langle x_{n}\rangle$ of any hypernode
in a galaxy other than $\Gamma_{0}$ must lie outside any finitely 
dispersed subgraph of $G$ for almost all $n$ whatever be the choice of that 
finitely dispersed subgraph. 

For instance, consider

{\bf Example 3.6.}  Let $D_{2}$ be the 2-dimensional grid;  that is,
we can represent $D_{2}$ by having its nodes at the lattice points 
$(k,l)$ of the 2-dimensional plane, where $k,l\in \Z$ and with its 
branches being $\{(k,l),(k+1,l)\}$ and $\{(k,l),(k,l+1)\}$.
So, the hypernodes of $^{*}\!D_{2}$ occur at $^{*}\!\Z\,\times\,^{*}\!\Z$.
Under this representation, the principal nodal galaxy of 
$^{*}\!D_{2}$ will have its nodes at the lattice points of 
$\Z\times\Z$.

Next, let $G$ be a connected graph obtained from $D_{2}$ by deleting 
or appending finitely many branches to $D_{2}$.  So, 
outside a finitely dispersed subgraph of $G$, $G$ is identical
to $D_{2}$.  Then the principal galaxy $\Gamma_{0}$ of 
$^{*}\!G$ is the same as (i.e., is graphically isomorphic to) 
$G$, but every other galaxy is the same as $D_{2}$.  $\Box$

In view of Examples 3.3 and 3.4, the following theorem is pertinent.
As always, we assume that $G$ is connected and has an infinite node set $X$.

{\bf Theorem 3.7.}  {\em Let $G$ be locally finite.  Then, $^{*}\!G$
has at least one hypernode not in its principal galaxy $\Gamma_{0}$ and thus 
at least one galaxy $\Gamma_{1}$ different from $\Gamma_{0}$.}

{\bf Proof.}  Choose any $x_{0}\in X$.  By
connectedness and local finiteness,
for each $n\in\N$, the set $X_{n}$ of nodes that are at a distance of 
$n$ from $x_{0}$ is nonempty and finite.  Also, $\cup X_{n}\,=\,X$ 
by the connectedness of $G$.  By K\"{o}nig's Lemma \cite[page 40]{wi},
there is a one-ended path $P$ starting at $x_{0}$.  $P$ must
pass through every $X_{n}$. Thus, there is a subsequence 
$\langle x_{0},x_{1},x_{2},\ldots\rangle$ of the sequence of nodes of 
$P$ such that $x_{n}\in X_{n}$; that is, $d(x_{n},x_{0})=n$ for every $n$.
Set ${\bf x}=[x_{n}]$.  Then, ${\bf x}$ must be in a galaxy 
$\Gamma_{1}$ that is different from the principal galaxy $\Gamma_{0}$.
$\Box$

\section{When $^{*}\!G$ Has a Hypernode Not in Its Principal Galaxy}

In this section, $G$ is connected and infinite but not necessarily 
locally finite.  Let $\Gamma_{a}$ and $\Gamma_{b}$ be two galaxies that are 
different from the principal galaxy $\Gamma_{0}$ of $^{*}\!G$.
We shall say that $\Gamma_{a}$ {\em is closer to $\Gamma_{0}$
than is} $\Gamma_{b}$ and that $\Gamma_{b}$ {\em is further away from 
$\Gamma_{0}$ than is} $\Gamma_{a}$
if there are a ${\bf y}=[y_{n}]$ in $\Gamma_{a}$ and a 
${\bf z}=[z_{n}]$ in $\Gamma_{b}$ such that, 
for some ${\bf x}=[x_{n}]$ in $\Gamma_{0}$
and for every $m\in\N$, we have
\[ N_{0}(m)\;=\;\{n\!:d(z_{n},x_{n})-d(y_{n},x_{n})
\,\geq\, m\}\;\in\;{\cal F}. \]
Any set of galaxies for which every two of them, say,
$\Gamma_{a}$ and $\Gamma_{b}$  satisfy this condition will be said to be
{\em totally ordered according to their closeness to}
$\Gamma_{0}$.  With Lemma 3.1 in hand, the conditions for a total ordering 
(reflexivity, antisymmetry, transitivity, and 
connectedness) are readily shown.  For instance, the proof of Theorem 4.3
below establishes transitivity.

{\bf Lemma 4.1.}  {\em These definitions are independent of the 
representative sequences $\langle x_{n}\rangle$, $\langle y_{n}\rangle$, 
and $\langle z_{n}\rangle$ chosen for ${\bf x}$, ${\bf y}$,
and ${\bf z}$.}

{\bf Proof.} Let $\langle x_{n}'\rangle$, $\langle y_{n}'\rangle$, 
and $\langle z_{n}'\rangle$ be any other such representative sequences.
Then, 
\[ d(z_{n},x_{n})\;\leq\;d(z_{n},z_{n}')+d(z_{n}',x_{n}')+d(x_{n}',x_{n}). \]
So, 
\[ d(z_{n}',x_{n}')\;\geq\;d(z_{n},x_{n})-d(z_{n},z_{n}')-d(x_{n}',x_{n})
\;=\; d(z_{n},x_{n}) \]
for all $n$ in some $N_{1}\in{\cal F}$.
Also, 
\[ d(y_{n}',x_{n}')\;\leq\; 
d(y_{n}',y_{n})+d(y_{n},x_{n})+d(x_{n},x_{n}')\;=\; d(y_{n},x_{n}) \]
for all $n$ in some $N_{2}\in{\cal F}$.
Therefore, 
\[ d(z_{n}',x_{n}')-d(y_{n}',x_{n}')\;\geq\; d(z_{n},x_{n})-d(y_{n},x_{n}) \]
for all $n$ in $N_{1}\cap N_{2}\,\in\, {\cal F}$.
So, for $N_{0}(m)$ as defined above and for each $m$ 
no matter how large,
\[ \{n\!: d(z_{n}',x_{n}')-d(y_{n}',x_{n}')\;\geq\; m\}
\;\supseteq\; N_{0}(m)\cap N_{1}\cap N_{2}\;\in\; {\cal F}, \]
which implies that the left-hand side is also a set in $\cal F$.
This proves Lemma 4.1.  $\Box$

We will say that a set $A$ is a {\em totally ordered, two-way infinite 
sequence} if there is a bijection from the set $\Z$ of integers to the set 
$A$ that preserves the total ordering of $\Z$.

{\bf Theorem 4.2.}  {\em If $^{*}\!G$ has a hypernode 
${\bf v}=[v_{n}]$ that is not in 
its principal galaxy $\Gamma_{0}$, then there exists a two-way 
infinite sequence of galaxies totally ordered according to their 
closeness to $\Gamma_{0}$ and with $\bf v$ being in one of those galaxies.}

{\bf Note.} There may be many such sequences, and a galaxy in one 
sequence and a galaxy in another sequence 
may not be comparable according to their 
closeness to $\Gamma_{0}$. Also, a somewhat 
different version of this theorem with a rather longer proof 
can be found in the archival website, www.arxiv.org, 
under Mathematics, Zemanian.

{\bf Proof.}  Let ${\bf x}=[\langle x,x,x,\ldots\rangle ]$ be a 
standard hypernode in $\Gamma_{0}$.  Also, let ${\bf v}=[v_{n}]$
be the asserted hypernode not in $\Gamma_{0}$.  Thus, for each 
$m\in\N$, 
\begin{equation}
\{n\!: d(x,v_{n})>m\}\,\in\,{\cal F}.  \label{4.1}
\end{equation}
Between every two nodes of a connected, conventionally infinite 
graph there is a geodesic path whose length is equal to the 
distance between those nodes.  For each $n\in\N$, choose a 
geodesic path $P$ in $G$ terminating at $x$ and $v_{n}$.
If the natural  number $d(x,v_{n})$ is even (resp. odd), there is a 
unique node $u_{n}$ in $P$ such that $d(x,u_{n})=d(u_{n},v_{n})
=d(x,v_{n})/2$ (resp. $d(x,u_{n})=d(u_{n},v_{n})-1=(d(x,v_{n})-1)/2)$.
It follows from this that, if there is a $k\in\N$ such that 
$\{n\!: d(x,u_{n})\leq k\}\in{\cal F}$, then there is a 
$k'\in\N$ with $\{n\!: d(x,v_{n})\leq k'\}\in{\cal F}$,
in violation of (\ref{4.1}). Consequently, for each 
$m\in\N$, $\{n\!: d(x,u_{n})>m\}\in{\cal F}$.  This implies that 
${\bf u}=[u_{n}]$ is in a galaxy $\Gamma_{a}$ different
from the principal galaxy $\Gamma_{0}$.

Furthermore, with $d(x,v_{n})$ being even (resp. odd) again and
with $u_{n}$ being chosen as before, 
\[ d(x,v_{n})\,-\, d(x,u_{n})\;=\; d(x,v_{n})/2 \]
(resp.
\[ d(x,v_{n})\,-\,d(x,u_{n})\;=\;(d(x,v_{n})+1)/2). \]
Now, if there is a $k\in\N$ such that 
\[ \{n\!: d(x,v_{n})-d(x,u_{n})\leq k\}\;\in\;{\cal F}, \]
then there is a $k'\in\N$ such that 
\[ \{n\!: d(x,v_{n})\leq k'\}\;\in\;{\cal F}, \]
again in violation of (\ref{4.1}).  Thus, for each $m\in\N$,
\[ \{n\!: d(x,v_{n})\,-\,d(x,u_{n})\,>\,m\}\;\in\;{\cal F}.\]
This implies that ${\bf u}=[u_{n}]$ and ${\bf v}=[v_{n}]$
are in different galaxies, $\Gamma_{a}$ and $\Gamma_{b}$
respectively, with $\Gamma_{a}$ being closer to $\Gamma_{0}$ than
is $\Gamma_{b}$.  

We can now repeat this argument with $\Gamma_{b}$ replaced
by $\Gamma_{a}$ and with ${\bf u}=[u_{n}]$ playing the role that
${\bf v}=[v_{n}]$ played
to find still another galaxy $\Gamma_{a}'$ different from 
$\Gamma_{0}$ and closer to $\Gamma_{0}$ than is $\Gamma_{a}$.  
Continual repetitions yield an infinite sequence of galaxies 
indexed by, say, the negative integers and totally ordered
by their closeness to $\Gamma_{0}$.

The conclusion that there is an infinite sequence of galaxies
progressively further away from $\Gamma_{0}$ than is 
$\Gamma_{b}$ is easier to prove.  With ${\bf v}\in \Gamma_{b}$ as before,
we have that, for every $m\in\N$, $\{n\!: d(x,v_{n})>m\}\in {\cal F}$.
Therefore, for each $n\in\N$, we can choose $w_{n}$ as an element of 
$\langle v_{n}\rangle$ such that 
\begin{equation}
d(x,w_{n})\geq d(x,v_{n})+n. \label{4.2}
\end{equation}
Hence, for each $m\in\N$, 
\[ \{n\!: d(x,w_{n})>m\}\;\supseteq\; \{n\!: d(x, v_{n}) > m\}. \]
Since the right-hand side is a member of $\cal F$, 
so too is the left-hand side.  Thus, ${\bf w}=[w_{n}]$ is in a 
galaxy different from $\Gamma_{0}$.  Moreover,
from (\ref{4.2}) we have that, for each $m\in\N$,
$\{n\!: d(x,w_{n})-d(x,v_{n})>m\}$ is a cofinite set and therefore 
is a member of $\cal F$.
Consequently, $\bf w$ is in a galaxy $\Gamma_{c}$ that is further away from
$\Gamma_{0}$ than is $\Gamma_{b}$.

We can repeat the argument of the last paragraph with
$\Gamma_{c}$ in place of $\Gamma_{b}$ to find still another galaxy
$\Gamma_{c}'$ further away from $\Gamma_{0}$ than is $\Gamma_{c}$.
Repetitions of this argument show that there is an 
infinite sequence of galaxies indexed by, say, the 
positive integers and totally 
ordered by their closeness to $\Gamma_{0}$.
  
The conjunction of the two infinite sequences along with $\Gamma_{b}$
yields the conclusion of 
the theorem.  $\Box$

By virtue of Theorem 3.7, the 
conclusion of Theorem 4.2 holds whenever $G$ is locally finite.

In general, the hypothesis of Theorem 4.2 may or may not hold.  Thus, 
$^{*}\!G$ either has exactly one galaxy, its principal one $\Gamma_{0}$, 
or has infinitely many galaxies.

Instead of the idea of ``totally ordered according to closeness 
to $\Gamma_{0}$,'' we can define the idea of ``partially ordered 
according to closeness to $\Gamma_{0}$'' in much the same way.  Just drop the 
connectedness axiom for a total ordering.

{\bf Theorem 4.3.}  {\em Under the hypothesis of Theorem 4.2, 
the set of galaxies of $^{*}\!G$ is 
partially ordered according to the closeness of the galaxies 
to the principal galaxy $\Gamma_{0}$.}

{\bf Proof.}  Reflexivity and antisymmetry are obvious.  
Consider transitivity:  Let $\Gamma_{a}$, $\Gamma_{b}$, and $\Gamma_{c}$ 
be galaxies different from $\Gamma_{0}$. (The case where 
$\Gamma_{a}=\Gamma_{0}$ can be argued similarly.)
Assume that $\Gamma_{a}$ is closer to $\Gamma_{0}$ than is 
$\Gamma_{b}$ and that $\Gamma_{b}$ is closer to $\Gamma_{0}$
than is $\Gamma_{c}$.  Thus, for any ${\bf x}$ in $\Gamma_{0}$,
${\bf u}$ in $\Gamma_{a}$, ${\bf v}$ in $\Gamma_{b}$, and 
${\bf w}$ in $\Gamma_{c}$ and for every $m\in\N$, we have
\[ N_{uv}\;=\;\{n\!: d(v_{n},x_{n})-d(u_{n},x_{n})\geq m\}\,\in\,{\cal F} \]
and 
\[ N_{vw}\;=\;\{n\!: d(w_{n},x_{n})-d(v_{n},x_{n})\geq m\}\,\in\,{\cal F}. \]
We also have
\[ d(w_{n},x_{n})-d(u_{n},x_{n})\;=\;d(w_{n},x_{n})-d(v_{n},x_{n})
+d(v_{n},x_{n})-d(u_{n},x_{n}). \]
So, 
\[ N_{uw}\;=\;\{n\!: d(w_{n},x_{n})-d(u_{n},x_{n})\geq 2m\}\,
\supseteq\,N_{uv}\cap N_{vw}\;\in\;{\cal F}. \]
Thus, $N_{uw}\in {\cal F}$. Since $m$ can be chosen arbitrarily,
we can conclude that $\Gamma_{a}$ is closer to $\Gamma_{0}$ 
than is $\Gamma_{c}$.  $\Box$

\section{The Hyperordinals}

In the following sections, we shall extend the results obtained so far 
to enlargements of transfinite graphs of rank 1, 
that is, to enlargements of 1-graphs.  For this purpose, 
we need to replace the set $^{*}\!\N$ of hypernaturals by
a set of ``hyperordinals'';  these are defined as follows.
A hyperordinal $\underline{\alpha}$ is an equivalence class 
of sequences of ordinals where two such sequences 
$\langle \alpha_{n}\rangle$ and $\langle \beta_{n}\rangle$ are taken
to be equivalent if $\{n\!: \alpha_{n}=\beta_{n}\}\in {\cal F}$.
We denote $\underline{\alpha}$ also by $[\alpha_{n}]$ 
where again the $\alpha_{n}$
are the elements of one (any one) of the sequences in the equivalent class.
Any set of hyperordinals is totally ordered by the inequality relation.
That is, given any hyperordinals $\underline{\alpha}=[\alpha_{n}]$ 
and $\underline{\beta}=[\beta_{n}]$,
exactly one of the sets:
\[ \{n\!: \alpha_{n}<\beta_{n}\},\;\; \{n\!:\alpha_{n}=\beta_{n}\},\;\; \{n\!:\alpha_{n}>\beta_{n}\} \]
will be in ${\cal F}$.  So, exactly one of the expressions:
\[ \underline{\alpha}<\underline{\beta},\;\;\underline{\alpha}=\underline{\beta},\;\;\underline{\alpha}>\underline{\beta} \] 
holds.

\section{Walks in 1-Graphs}

1-graphs arise when conventionally infinite graphs 
are connected at their infinite extremities 
through 1-nodes, the latter being a generalization of the
idea of a node.  Such 1-nodes and the resulting 1-graphs are defined in 
\cite[Section 2.1]{tgen} and also in \cite[Section 2.3]{gn}.  
Let us restate the needed definitions
concisely.

We will be dealing with two kinds of nodes and two kinds of graphs.
A conventionally infinite graph $G^{0}$
will now be called a 0-{\em graph} and the nodes in
$G^{0}$ will be called 0-{\em nodes} in order to distinguish these ideas 
from those pertaining to transfinite graphs of rank 1.
Similarly, what we called a ``hypernode'' previously will henceforth be called 
a 0-{\em hypernode}, and what we called a ``galaxy'' in the enlargement of
a 0-graph will now be called a 0-{\em galaxy}. 

An {\em infinite extremity} of a 0-graph $G^{0}$ is defined as an 
equivalence class of one-ended paths in $G^{0}$, where two such paths
are considered to be {\em equivalent} if they are eventually identical.
Such an equivalence class is called a 0-{\em tip} of $G^{0}$.
$G^{0}$ may have one or more 0-tips (or possibly none at all).
To obtain the ``1-nodes,'' the set of 0-tips is partitioned 
in some fashion into subsets, and to each subset a single 0-node
may (or may not) be added under the proviso that, if a 0-node
is added to one subset, it is not added to any other subset.
Then, each subset (possibly augmented with a 0-node) is called 
a 1-{\em node}.  With $X^{1}$ denoting the set of 1-nodes and $X^{0}$
the set of 0-nodes of $G^{0}$, the 1-{\em graph}  $G^{1}$
is defined as the triplet:
\[G^{1}\;=\;\{X^{0},B,X^{1}\}, \]
and $G^{0}=\{X^{0},B\}$ is now called the 0-{\em graph of}
$G^{1}$.  Furthermore, a path in $G^{0}$ is now called a 0-{\em path}, 
and connectedness in $G^{0}$ is now called 0-{\em connectedness}.
We will consistently append the superscript 0 to the symbols and 
the prefix 0- to the 
terminology for concepts from Sections 2 through 4 regarding 0-graphs.

In order to define the ``1-galaxies,'' we need the idea
of distances in a 1-graph $G^{1}$.  But now, we must make 
a significant choice.  The  distances between two nodes 
(0-nodes or 1-nodes) can be defined as the minimum length of 
all paths---or, alternatively, of all walks---connecting the two nodes.  
It turns out that a path need not exist between two nodes
in a 1-graph $G^{1}$, but a walk always will exist between them. 
To ensure the existence of at least one path between every two nodes, 
additional conditions must be imposed on $G^{1}$ (see 
\cite[Conditions 3.2-1 and 3.5-1]{tgen} or \cite[Condition 3.1-2]{gn}), 
and this leads to a more restrictive and yet more complicated 
theory involving distances.  Such can be done,
but it is more general and simpler to use walk-based distance ideas.
This we now do.

A {\em nontrivial 0-walk} $W^{0}$ in a 0-graph is the conventional concept.
It is a (finite or one-way infinite or two-way infinite) 
alternating sequence:
\begin{equation}
W^{0}\;=\;\langle \ldots,x_{-1}^{0},b_{-1},x_{0}^{0},b_{0},x_{1}^{0},b_{1},\ldots\rangle  \label{6.1}
\end{equation}
of 0-nodes $x_{m}^{0}$ and branches $b_{m}$,
where each branch $b_{m}$ is incident to the two 0-nodes $x_{m}^{0}$ 
and $x_{m+1}^{0}$ adjacent to it in the sequence.  
If the sequence terminates at either side, it is required to terminate
at a 0-node. The 0-walk is called {\em two-ended} or 
{\em finite} if it terminates on both sides, {\em one-ended} if it 
terminates on just one side, and {\em endless} 
if it terminates on neither side.

A {\em trivial 0-walk} is a singleton set whose sole element is a 0-node.

A one-ended 0-walk $W^{0}$ will be called {\em extended} if its 0-nodes 
are eventually distinct, that is, if it is eventually identical
to a one-ended path.  We say that $W^{0}$ {\em traverses} a 0-tip
if it is extended and eventually identical to
a representative of that 0-tip.  Finally, $W^{0}$ is said
to {\em reach} a 1-node $x^{1}$ if $W^{0}$ traverses 
a 0-tip contained in $x^{1}$.  In the same way,
an endless 0-walk can {\em reach} two 1-nodes
(or possibly reach the same 1-node) by traversing two 
0-tips, one toward the left and the other toward the right.
When this is so, we say that the endless 0-walk is {\em extended}.
On the other hand, if a 0-walk terminates at a 0-node contained in
a 1-node, we again say that the 0-walk {\em reaches}
both of those nodes and does so {\em through} a branch incident to that 
0-node.

Every two-ended 0-walk contains a 0-path that terminates at the two 
0-nodes at which the 0-walk terminates, so there is no need to employ 
0-walks when defining distances in a 0-graph.  On the other hand, 
such a need arises for 1-graphs.  To meet this need, we 
first define a 0-{\em section} $S^{0}$ in a 1-graph $G^{1}$
as a subgraph $S^{0}$ of the 0-graph $G^{0}$ of $G^{1}$
induced by a maximal set of branches that are pairwise 0-connected in 
$G^{0}$.  A 1-node $x^{1}$ is said to be {\em incident to} 
$S^{0}$ if either it contains a 0-node incident to a branch
of $S^{0}$ or it contains a 0-tip having a representative one-ended path 
lying entirely within $S^{0}$.  In this case, we also say
that that 0-tip {\em belongs to} $S^{0}$.
Given two 1-nodes $x^{1}$ and $y^{1}$ incident to $S^{0}$, there 
will be a 0-walk $W^{0}$ in $S^{0}$ that reaches each of 
$x^{1}$ and $y^{1}$ through a 0-tip belonging to $S^{0}$ or 
through a branch in $S^{0}$.\footnote{For examples of when a 0-walk is needed
because a 0-path won't do, see Figures 3.1 and 3.2 of \cite{tgen}
and Figures 4.1, 5.1, 5.2, and 5.3 of \cite{gn}.}
Moreover, there may also be a 0-walk $W^{0}$ in $S^{0}$ that reaches the 
same 1-node at both extremities of $W^{0}$. To be more specific, let 
us state

{\bf Lemma 6.1.}  {\em Let $S^{0}$ be a 0-section in $G^{1}$, and let
$x^{1}$ and $y^{1}$ be two 1-nodes incident to $S^{0}$.  
Then, there exists a 0-walk in $S^{0}$ that reaches 
$x^{1}$ and $y^{1}$.}

{\bf Proof.}  That $x^{1}$ is incident to $S^{0}$ means that
there is a 0-path $P_{x}^{0}$ in $S^{0}$ that either reaches $x^{1}$
through a 0-tip of $x^{1}$ or reaches $x^{1}$ through a branch.
Similarly, there is such a 0-path $P_{y}^{0}$ reaching $y^{1}$.  Let 
$u^{0}$ be a node of $P_{x}^{0}$, and let $v^{0}$ be a node 
of $P_{y}^{0}$.  Since $S^{0}$ is 0-connected, 
there is a 0-path $P_{uv}^{0}$ in $S^{0}$ terminating at $u^{0}$ and $v^{0}$
(possibly a trivial 0-path if $u^{0}=v^{0}$).
Then, $P_{x}^{0}\cup P_{uv}^{0}\cup P_{y}^{0}$ as a 0-walk 
in $S^{0}$ as asserted. $\Box$

A {\em nontrivial, two-ended 1-walk} $W^{1}$ is a finite sequence:
\begin{equation}
W^{1}\;=\;\langle x_{0},W^{0}_{0},x_{1}^{1},W_{1}^{0},\ldots,x_{m-1}^{1},W_{m-1}^{0},x_{m}\rangle  \label{6.2} 
\end{equation}
with $m\geq 1$ that satisfies the following conditions.
\begin{description}
\item{1.} $x_{1}^{1},\ldots,x_{m-1}^{1}$ are 1-nodes, while $x_{0}$ 
and $x_{m}$ may be either 0-nodes or 1-nodes.
\item{2.} For each $k=0,\ldots,m-1$, $W_{k}^{0}$ is a nontrivial 0-walk 
that reaches the two nodes adjacent to it in the sequence.
\item{3.} For each $k=1,\ldots,m-1$, at least one of 
$W_{k-1}^{0}$ and $W_{k}^{0}$ reaches $x_{k}^{1}$
through a 0-tip, not through a branch.  Also, if $m=1$, $W_{0}^{0}$
reaches at least one of $x_{0}$ and $x_{1}$ through a 0-tip.  
\end{description}
A {\em one-ended} 1-walk is a sequence like (\ref{6.2})
except that it extends infinitely to the right.  An {\em endless} 1-walk
extends infinitely on both sides. 
A {\em trivial 1-walk} is a singleton set whose sole element
is either a 0-node or a 1-node.

We now define a more general kind of connectedness (called ``1-wconnectedness''
to distinguish it from path-based 1-connectedness).
Two branches (resp. two nodes---either 0-nodes or 1-nodes)
will be said to be 1-{\em wconnected} if there exists a 0-walk 
or 1-walk that terminates at a 0-node of each branch (resp. that
terminates at those two nodes).  If a terminal node of a walk is 
the same as, or contains, or is contained in the terminal node of 
another walk, the two walks taken together 
form another walk.  We call this the {\em conjunction}
of the two walks.
It follows that 1-wconnectedness
is a transitive binary relation for the branch set $B$
of the 1-graph $G^{1}$ and is in fact an equivalence relation.
If every two branches of $G^{1}$ are 1-wconnected, 
we will say that $G^{1}$ is 1-{\em wconnected}.

\section{Walk-Based Distances in a 1-Graph}

The length $|W^{0}|$ of a 0-walk $W^{0}$ is defined as follows:
If $W^{0}$ is two-ended, $|W^{0}|$ is the number $\tau_{0}$
of branch traversals in it;  that is, each branch is counted as 
many times as it appears in $W^{0}$.
If $W^{0}$ is one-ended and extended, we set $|W^{0}|=\omega$, the
first transfinite ordinal.  If $W^{0}$ is endless and extended in 
both directions, we set $|W^{0}|=\omega\cdot 2$.

As for a nontrivial 
two-ended 1-walk $W^{1}$, its length $|W^{1}|$ is taken to be 
$|W^{1}|=\sum_{k=0}^{m}|W_{k}^{0}|$, where the sum is 
over the finitely many 0-walks $W_{k}^{0}$ in (\ref{6.2}).  Thus,
\begin{equation}
|W^{1}|\;=\;\omega\cdot \tau_{1}+\tau_{0}  \label{7.1}
\end{equation}
where $\tau_{1}$ is the number of traversals of 0-tips performed by $W^{1}$ 
and $\tau_{0}$ is the number of traversals of branches in all the 
two-ended (i.e., finite) 0-walks appearing as terms in (\ref{6.2}).
We take $\sum_{k=0}^{m} |W^{0}_{k}|$ to be the natural sum of ordinals; 
this yields a normal expansion of an ordinal \cite[pages 354-355]{ab}.
$\tau_{1}$ is not 0 because $W^{1}$ is a nontrivial, two-sided 1-walk.
However, $\tau_{0}$ may be 0, this occurring when every $W_{k}^{0}$
in (\ref{6.2}) is one-ended or endless.  

A 0-node is called {\em maximal} if it is not contained in a 1-node, 
and {\em nonmaximal} otherwise.  A distance measured from a 
nonmaximal 0-node is the same as that measured from the 1-node containing it.
Given two nodes $x$ and $y$ (of ranks 0 or 1), we define
the {\em wdistance}\footnote{We write ``wdistance'' to distinguish this 
walk-based idea from a distance based on paths.} $d(x,y)$ between 
them as 
\begin{equation}
d(x,y)\;=\;\min|W_{x,y}|  \label{7.2}
\end{equation}
where the minimum is taken over all two-ended walks  (0-walks or 1-walks)
terminating at $x$ and $y$.  That minimum exists because any set of 
ordinals is a well-ordered set.  In view of (\ref{7.1}), $d(x,y)<\omega^{2}$.
If $x=y$, we set $d(x,x)=0$. 

Clearly, if $x\neq y$, $d(x,y)>0$ and $d(x,y)=d(y,x)$.
Furthermore, the conjunction of two two-ended 
walks is again a two-ended walk, whose length is the natural sum 
of the ordinal lengths
of the two walks.  So, by taking minimums appropriately, we obtain the 
triangle inequality:
\begin{equation}
d(x,z)\;\leq\;d(x,y)\,+\,d(y,z)  \label{7.3}
\end{equation}
where again the natural sum of ordinals is understood.  
Altogether then, we have 

{\bf Lemma 7.1.} {\em The ordinal-valued wdistances between the maximal
nodes of a 1-graph satisfy the metric axioms.}

\section{Enlargements of 1-Graphs and Hyperdistances in Them}

In \cite[pages 163-164]{gn}, a nonstandard 1-node was defined 
as an equivalence class of sequences of sets of tips shorted together,
with the tips taken from sequences of possibly differing 1-graphs.
But, since each set of tips shorted together is a 1-node, 
that definition of a nonstandard 1-node can also be stated as 
an equivalence class of sequences of 1-nodes. 
Specializing to the case where all
the 1-graphs are the same, we have the following definition of a 
nonstandard 1-node, which we now call a ``1-hypernode.''

Consider a given 1-graph 
along with a chosen free ultrafilter ${\cal F}$.
Two sequences $\langle x_{n}^{1}\rangle$ and 
$\langle y_{n}^{1}\rangle$ of 1-nodes in $G^{1}$ are taken to be
{\em equivalent} if $\{n\!: x_{n}^{1}=y_{n}^{1}\}\in{\cal F}$.
It is easy to show that this is truly an equivalence relation.
Then, ${\bf x}^{1}=[x_{n}^{1}]$ denotes one such equivalence class,
where the $x_{n}^{1}$ are the elements of any one of the
sequences in that class.  ${\bf x}^{1}$ will be called 
a 1-{\em hypernode}.

The {\em enlargement} of the 1-graph $G^{1}=\{X^{0},B,X^{1}\}$
is the nonstandard 1-graph
\[ ^{*}\!G^{1}\;=\;\{\,^{*}\!X^{0},\,^{*}\!B,\,^{*}\!X^{1}\,\} \]
where $^{*}\!X^{0}$ and $^{*}\!B$ are respectively the set of 
0-hypernodes and branches in the enlargement of the 0-graph
$G^{0}=\{X^{0},B\}$ of $G^{1}$ and $^{*}\!X^{1}$
is the set of 1-hypernodes defined above, that is, the set of all 
equivalence classes of sequences of 1-nodes taken from $X^{1}$.

We define the {\em hyperdistance} ${\bf d}$ between any two hypernodes
${\bf x}$ and ${\bf y}$ of $^{*}\!G^{1}$ (of ranks 0 and/or 1) 
to be the internal function
\begin{equation}
{\bf d}({\bf x},{\bf y})\;=\;[d(x_{n},y_{n})].  \label{8.1}
\end{equation}
Since distances in $G^{1}$ are less than $\omega^{2}$, 
${\bf d}({\bf x},{\bf y})$ is a hyperordinal 
less than $\underline{\omega}^{2}$.
We say that a 0-hypernode ${\bf x}^{0}=[x_{n}^{0}]$ is {\em maximal}
if the set of $n$ for which $x_{n}^{0}$ is not contained 
in a 1-node is a member of $\cal F$.  
All the 1-nodes in this work are perforce 
maximal because there are no nodes of higher rank.
${\bf d}$, when restricted to the maximal hypernodes,  
also satisfies the metric axioms, in particular, the 
triangle inequality:
\begin{equation}
{\bf d}({\bf x},{\bf z})\;\leq\;{\bf d}({\bf x},{\bf y})\,+\,{\bf d}({\bf y},{\bf z})  \label{8.2}
\end{equation}
But, now $\bf d$ is hyperordinal-valued.

\section{The Galaxies of $^{*}\!G^{1}$}

The 0-{\em galaxies} of $^{*}\!G^{1}$ are defined just as they are for
the enlargement $^{*}\!G^{0}$ of a 0-graph;  see Section 3.
However, we henceforth write ``0-galaxy'' in place of ``galaxy'' 
and ``0-limitedly distant'' in place of ``limitedly distant.''

As was mentioned above, each 0-section of $G^{1}$ is the subgraph 
of the 0-graph $G^{0}\{X^{0},B\}$ of $G^{1}$ induced by a 
maximal set of branches that are 0-connected.  
A 0-section is a 0-graph by itself.
So, within the 
enlargement $^{*}\!G^{1}$, each 0-section $S^{0}$ enlarges into 
$^{*}\!S^{0}$ as defined in Section 2.  Within each enlarged 
0-section there may be one or more 0-galaxies.  As a special case,
a particular 0-section may have only finitely many 
0-nodes, and so its enlargement is itself---all its 0-hypernodes 
are standard.  On the other hand, there may be infinitely
many 0-galaxies in some enlarged 0-section.  Moreover,
the enlarged 0-sections do not, in general, comprise all of the 
enlarged 0-graph $^{*}\!G^{0}=\{\,^{*}\!X^{0},\,^{*}\!B\,\}$
of $^{*}\!G^{1}$.  Indeed, there can be a 0-hypernode 
${\bf x}^{0}=[x_{n}^{0}]$ where each $x_{n}^{0}$ resides in a different
0-section; in this case ${\bf x}^{0}$ will reside in a 0-galaxy 
that is not in an enlargement of a 0-section.  

Something more can happen with regard to the 0-galaxies in 
$^{*}\!G^{1}$.  0-galaxies can now contain 1-hypernodes.  
For example, this occurs when a 1-node $x^{1}$ is incident
to a 0-section $S^{0}$ through a branch.  Then, the standard 1-hypernode 
${\bf x}^{1}$ corresponding to $x^{1}$ is 0-limitedly distant from the 
standard 0-hypernodes in $^{*}\! S^{0}$.  So, 
there is a 0-galaxy containing
not only $^{*}\!S^{0}$ but ${\bf x}^{1}$ as well.  See Example 9.3 below in
this regard.  In general, the nodal 0-galaxies partition the set 
$^{*}\!X^{0}\,\cup\,^{*}\!X^{1}$ of all the hypernodes in 
$^{*}\!G^{1}$.  As we shall see in Examples 9.1 and 9.2 below, there
may be a singleton 0-galaxy containing a 1-hypernode only. 

Let us now turn to the ``1-galaxies'' of $^{*}\!G^{1}$
Two hypernodes ${\bf x}=[x_{n}]$ and 
${\bf y}=[y_{n}]$ (of ranks 0 and/or 1) in 
$^{*}\!G^{1}$ will be said to be in the same {\em nodal 1-galaxy}
$\dot{\Gamma}^{1}$ if there exists a natural number $k\in\N$
depending on the choices of $\bf x$ and $\bf y$
such that $\{n\!:d(x_{n},y_{n})\leq \omega\cdot k\}\in{\cal F}$.
In this case, we say that ${\bf x}$ and ${\bf y}$ are 
{\em 1-limitedly distant}, and we write ${\bf d}({\bf x},{\bf y})
\leq [\omega\cdot k]$ where $[\omega\cdot k]$ denotes the standard 
hyperordinal corresponding to $\omega\cdot k$.
This defines an equivalence relation on the set 
$^{*}\!X^{0}\,\cup\, ^{*}\!X^{1}$ of all the hypernodes in $^{*}\!G^{1}$.
Indeed, reflexivity and symmetry are obvious.  For transitivity, 
assume that ${\bf x}$ and ${\bf y}$ are 1-limitedly distant and that
${\bf y}$ and ${\bf z}$ are 1-limitedly distant, too.  Then,
there are natural numbers $k_{1}$ and $k_{2}$ such that 
\[ N_{xy}\;=\;\{n\!: d(x_{n},y_{n})\leq \omega\cdot k_{1}\}\,\in\,{\cal F} \]
and 
\[ N_{yz}\;=\;\{n\!: d(y_{n},z_{n})\leq \omega\cdot k_{2}\}\,\in\,{\cal F}. \]
By the triangle inequality (\ref{7.3}),
\[ N_{xz}\;=\;\{n\!:d(x_{n},z_{n})\leq\omega\cdot (k_{1}+k_{2})\}\;\supseteq
\;N_{xy}\cap N_{yz}\;\in\;{\cal F}. \]
So, $N_{xz}\in {\cal F}$ and therefore ${\bf x}$ and ${\bf z}$ are
1-limitedly distant.  We can conclude that the set 
$^{*}\!X^{0}\,\cup\,^{*}\!X^{1}$
of all hypernodes in $^{*}\! G^{1}$
is partitioned into nodal 1-galaxies by this equivalence relation.

Corresponding to each nodal 1-galaxy $\dot{\Gamma}^{1}$, we
define a 1-{\em galaxy} $\Gamma^{1}$ as a nonstandard subgraph of 
$^{*}\!G^{1}$ consisting of all the hypernodes in 
$\dot{\Gamma}^{1}$ along with all the hyperbranches both of whose 
0-hypernodes are in $\dot{\Gamma}^{1}$.

No hyperbranch can have its two incident 0-hypernodes in two different 
0-galaxies or two different 1-galaxies because the distance between their 
0-hypernodes is 1.  Thus, the hyperbranch set $^{*}\!B$ is also partitioned 
by the 0-galaxies and more coarsely
by the 1-galaxies.  

The {\em principal 1-galaxy} $\Gamma_{0}^{1}$ of $^{*}\!G^{1}$
is the 1-galaxy whose hypernodes are 1-limitedly distant from a standard
hypernode in $^{*}\!G^{1}$ (i.e., from a node of $G^{1}$).

Note that the enlargement $^{*}\!S^{0}$ of each 0-section
$S^{0}$ of $G^{1}$ has its own principal 0-galaxy 
$\Gamma_{0}^{0}(S^{0})$.  Moreover, every $^{*}\!S^{0}$ lies 
within the principal 1-galaxy $\Gamma_{0}^{1}$. 
Indeed, any standard hypernode ${\bf x}$ by which
$\Gamma_{0}^{1}$ may be defined and any standard 0-hypernode ${\bf y}^{0}$
by which $\Gamma_{0}^{0}(S^{0})$ may be defined are 1-limitedly 
distant.  Also, the hyperdistance ${\bf d}({\bf y}^{0},{\bf z}^{0})$
between any two 0-hypernodes ${\bf y}^{0}$ and ${\bf z}^{0}$ of 
$^{*}\!S^{0}$ is no larger than a hypernatural ${\bf k}$.
So, by the triangle inequality (\ref{8.2}), every 0-hypernode of 
$^{*}\!S^{0}$ is 1-limitedly distant from ${\bf x}$.  Whence our assertion.

{\bf Example 9.1.}  Consider an endless 1-path $P^{1}$ having an
endless 0-path between every consecutive pair of 1-nodes in $P^{1}$.
The 0-sections of $P^{1}$ are those endless 0-paths, and each of their 
enlargements have an infinity of 
0-galaxies in $^{*}\!P^{1}$.  However, there are other
0-galaxies in $^{*}\!P^{1}$, infinitely many of them.  Indeed,
consider a 0-hypernode ${\bf x}^{0}=[x_{n}^{0}]$,
where each 0-node $x_{n}^{0}$ lies in a different
0-section of $P^{1}$;  ${\bf x}^{0}$ will lie in a 0-galaxy $\Gamma_{1}^{0}$
different from all the 0-galaxies in any
enlargement of a 0-section of $P^{1}$.  The 0-hypernodes of 
$\Gamma_{1}^{0}$ will be the 0-hypernodes that are 0-limitedly
distant from ${\bf x}^{0}$.  Furthermore, there are still other
0-galaxies now.  Each 1-hypernode ${\bf x}^{1}=[x_{n}^{1}]$ is the 
sole member of a 0-galaxy.  In fact, the nodal 0-galaxies partition 
the set of all the 0-hypernodes and 1-hypernodes.

On the other hand, the principal 1-galaxy of $^{*}\!P^{1}$ consists of all
the standard 1-hypernodes corresponding to the 
1-nodes of $P^{1}$ along with the enlargements of the 
0-sections of $P^{1}$.
Also, there will be infinitely many 1-galaxies, each of which
contains infinitely many 0-galaxies along with 1-hypernodes.  In
this particular case, each of the 1-galaxies is graphically
isomorphic to the principal 1-galaxy, but this is not true in general.  $\Box$

{\bf Example 9.2.}  An example of a nonstandard 1-graph 
$^{*}\!G^{1}$ having exactly one 1-galaxy 
(its principal one) and infinitely many 0-galaxies is 
provided by the enlargement of the 1-graph $G^{1}$
obtained from the 0-graph of 
Figure 1 by replacing each branch by an endless 0-path, 
thereby converting each 0-node into a 1-node.  Again each endless path
of that 1-graph $G^{1}$ is a 0-section, and its enlargement is like that 
of Example 3.2.  There are 
infinitely many such 0-galaxies in each enlargement of a 0-section.
Also, there are infinitely many 0-galaxies,
each consisting of a single 1-hypernode.
With regard to the 1-galaxies, the enlargement $^{*}\!G^{1}$
of $G^{1}$ 
mimics that of Example 3.4, except that now
the rank 0 is replaced by the rank 1. The hyperdistance between
every two 1-hypernodes (resp. 0-hypernodes) is no larger than
$\omega\cdot 4$ (resp. $\omega\cdot 6$).  Hence, $^{*}\!G^{1}$ has only 
one 1-galaxy, its principal one. $\Box$

{\bf Example 9.3.} Here is an example where each of the 1-hypernodes
is not isolated as the sole member
within a 0-galaxy.  Replace each of the 
horizontal branches in Figure 1 by an endless 0-path, but do not alter the
branches incident to $x_{g}$. Now, the nodes $x_{k}$ $(k=0,1,2,\ldots)$ 
become 1-nodes $x_{k}^{1}$, each containing a 0-node of the branch incident to 
$x_{k}^{1}$ and $x_{g}$.  Each 1-hypernode is 0-limitedly distant from the 
standard 0-hypernode ${\bf x}_{g}$ and thus is not so isolated.
The 1-hypernodes
(whether standard or nonstandard) along 
with the standard 0-hypernode for $x_{g}$ and the 
(standard or nonstandard) hyperbranches 
incident to ${\bf x}_{g}$ 
all comprise a single 0-galaxy.  Moreover, there will 
be other 0-galaxies obtained through equivalence classes of sequences of 
these nodes and branches.  In fact, each of the 
endless paths that replace the horizontal
branches yield infinitely many 0-galaxies.  Again, the nodal 
0-galaxies partition the set of all the hypernodes in $^{*}\!G^{1}$.

On the other hand, there is again only one 
1-galaxy for $^{*}\!G^{1}$.  $\Box$

{\bf Example 9.4.}  The distances in the three preceding examples
can be fully defined by paths.
So, let us now present an example where walks are needed.  The 1-graph
$G^{1}$ of Figure 3 illustrates one such case.  It consists of an 
infinite sequence of 0-subgraphs, each of which is an  
infinite series connections of four-branch subgraphs, each
in a diamond configuration, as shown.  To save words,
we shall refer to such an infinite series connection as a ``chain.''
The chain starting at the 0-node $x_{k}^{0}$ will be denoted by $C_{k}$
$(k=0,1,2,\ldots)$.  Each $C_{k}$ is a 0-graph;  it does not contain 
any 1-node.  Each $C_{k}$ has uncountably many 0-tips.  One 0-tip
has a representative 0-path starting at $x_{k}^{0}$, 
proceeding along the left-hand sides of 
the diamond configurations, and reaching the 1-node $x_{k}^{1}$.  
Another 0-tip has a representative 0-path 
that proceeds along the right-hand sides
and reaches the 1-node $x_{k+1}^{1}$.  Still other 0-tips 
of $C_{k}$ (uncountably
many of them) have representatives that pass back and forth 
between the two sides infinitely often
to reach singleton 1-nodes; these are not shown
in that figure.  The chain $C_{k}$ is connected to $C_{k+1}$ through the 
1-node $x_{k+1}^{1}$, as shown. Note that there is no path connecting, say, 
$x_{k}^{0}$ to $x_{m}^{0}$ when $m-k\geq 2$, but there is such a walk.

Each $C_{k}$ is a 0-section, and its enlargement $^{*}\!C_{k}$
has infinitely many 0-galaxies.  Also, the 1-nodes 
$x_{k}^{1}$ together produce infinitely many 0-galaxies,
each being a single 1-hypernode.  As before, the nodal
0-galaxies comprise a partition of $^{*}\!X^{0}\,\cup\,^{*}\!X^{1}$.

On the other hand, the enlargement $^{*}\!G^{1}$ of the 1-graph
$G^{1}$ of Figure 3 has infinitely many 1-galaxies.  Its principal
one is a copy of $G^{1}$.  Each of the other 1-galaxies is 
also a copy of $G^{1}$ except that it extends infinitely in both
directions---infinitely to the left and infinitely to the right.
Here, too, the nodal 1-galaxies comprise a partitioning of 
$^{*}\!X^{0}\,\cup\,^{*}\!X^{1}$, but a coarser one.  $\Box$

These examples indicate that the enlargements of 1-graphs 
can have rather complicated structures.

\section{Locally 1-Finite 1-Graphs and a Property of Their Enlargements}

In general, $^{*}\!G^{1}$ has 1-galaxies other than its principal 1-galaxy.
One circumstance where this occurs is when $^{*}\!G^{1}$
is locally finite in certain way, which we will explicate 
below. 

We need some 
more definitions.  Two 1-nodes of $G^{1}$ are 
said to be 1-{\em wadjacent} if they are incident to the same 0-section.  
A 1-node will be called a {\em boundary 1-node} if it is incident to 
two or more 0-sections.  $G^{1}$ will be called {\em locally
1-finite} if each of its 0-sections has only finitely many incident 
boundary 1-nodes.\footnote{Note that a 0-section in a locally 1-finite 1-graph 
may have infinitely many incident 1-nodes that are not boundary 1-nodes.
Also, this definition of locally 1-finiteness does not prohibit 
0-nodes of infinite degree.}  

{\bf Lemma 10.1.}  {\em Let $x^{1}$ be a boundary 1-node.  Then,
any 1-walk that passes through $x^{1}$ from any 0-section $S^{0}_{1}$ 
incident to $x^{1}$ to any other 0-section $S^{0}_{2}$ 
incident to $x^{1}$ must have a length no less than $\omega$.}

{\bf Proof.}  The only way such a walk can have a length less than $\omega$
(i.e., a length equal to a natural number) is if it avoids traversing a 0-tip 
in $x^{1}$.  But, this means that it passes through two branches incident 
to a 0-node in $x^{1}$.  But, that in turn means that $S^{0}_{1}$ and 
$S^{0}_{2}$ cannot be different 0-sections.  $\Box$

Remember that $G^{1}$ is called 
1-wconnected if, for every two nodes of $G^{1}$, there is a 0-walk 
or 1-walk that reaches those two nodes.

{\bf Lemma 10.2.}  {\em Any two 1-nodes $x^{1}$ and $y^{1}$ 
that are 1-wconnected but are not 1-wadjacent must satisfy
$d(x^{1},y^{1})\geq \omega$.}

{\bf Proof.} Any walk 1-wconnecting $x^{1}$ and $y^{1}$ must
pass through at least one boundary 1-node different from $x^{1}$ and $y^{1}$ 
while passing from one 0-section to another 0-section. 
Therefore, that walk must be a 1-walk.  By Lemma 10.1,
its length is no less than $\omega$.  Since this is true 
for every such walk, our conclusion follows.  $\Box$

The next theorem mimics Theorem 3.7 but at the rank 1.

{\bf Theorem 10.3.}  {\em Let $G^{1}$ be locally 1-finite and 
1-wconnected and have infinitely many boundary 1-nodes.
Then, given any 1-node $x_{0}^{1}$ of $G^{1}$, there is
a one-ended 1-walk $W^{1}$ starting at $x_{0}^{1}$:
\[W^{1}\;=\;\langle x_{0}^{1},W_{0}^{0},x_{1}^{1},W_{1}^{0},
\ldots,x_{m}^{1},W_{m}^{0},\ldots\rangle \]
such that there is a subsequence of 1-nodes $x_{m_{k}}^{1}$, 
$k=1,2,3,\ldots$, satisfying 
$d(x_{0}^{1},x_{m_{k}}^{1})\,\geq\,\omega\cdot k$.}

{\bf Proof.}  $x^{1}_{0}$ need not be a boundary 1-node, but it will be 
1-wadjacent to only finitely many boundary 1-nodes because of 
local 1-finiteness and 1-wconnectedness.  Let $X_{0}$ be the 
nonempty finite set of those boundary 1-nodes.
For the same reasons, there is a nonempty finite 
set $X_{1}$ of boundary 1-nodes,
each being 1-wadjacent to some 1-node in $X_{0}$ but not 1-wadjacent 
to $x_{0}^{1}$.  By Lemma 10.2, for each $x^{1}\in X_{2}$, we have 
$d(x_{0}^{1},x^{1})\geq\omega$.  
In general, for each $k\in\N$, $k\geq 2$,
there is a nonempty finite set $X_{k}$ of boundary 1-nodes, each 
being 1-wadjacent to some 1-node in $X_{k-1}$ but not 1-wadjacent to any 
of the 1-nodes in $\cup_{l=0}^{k-2} X_{l}$.  
By Lemma 10.2 again, for any such $x^{1}\in X_{k}$, we have 
$d(x_{0}^{1},x^{1})\geq \omega\cdot k$.

We now adapt the proof of K\"{o}nig's lemma:  From each of the
infinitely any boundary 1-nodes in $G^{1}$, there is a 1-walk reaching that 
boundary 1-node and also reaching $x_{0}^{1}$.
Thus, there are infinitely many 1-walks starting at $x_{0}^{1}$
and passing through one of the 1-nodes in $X_{0}$, say, $x_{m_{0}}^{1}$.
Among those 1-walks, there are again infinitely many 1-walks
passing through one of the 1-nodes in $X_{1}$, say, $x_{m_{1}}^{1}$.
Continuing in this say, we find an infinite sequence 
$\langle x_{m_{1}}^{1}, x_{m_{2}}^{1},x_{m_{3}}^{1},\ldots\rangle$
of 1-nodes occurring in a one-ended 1-walk starting at 
$x_{0}^{1}$ and such that $d(x_{0}^{1},x_{m_{k}}^{1})\geq\omega
\cdot k$. $\Box$

{\bf Corollary 10.4.}  {\em Under the hypothesis of Theorem 10.3, the enlargement 
$^{*}\!G^{1}$ of $G^{1}$ has at least one 
1-hypernode not in its principal galaxy $\Gamma_{0}^{1}$ and 
thus at least one
1-galaxy $\Gamma^{1}$
different from its principal 1-galaxy $\Gamma_{0}^{1}$.}

{\bf Proof.}  Set ${\bf x}^{1} =[\langle x_{0}^{1},x_{m_{0}}^{1},
x_{m_{1}}^{1},\ldots\rangle]= [x_{m_{n}}^{1}]$, 
where the $x_{m_{n}}^{1}$ are the 1-nodes
specified in the preceding proof (replace $k$ by $n$).  
With ${\bf x}_{0}^{1}$
being the standard 1-hypernode corresponding to $x_{0}^{1}$, we have by 
Theorem 10.3 that ${\bf d}({\bf x}_{0}^{1},{\bf x}^{1})\geq 
[\omega\cdot n]$.  Hence, ${\bf x}^{1}$ is not 1-limitedly distant
from ${\bf x}_{0}^{1}$ and thus must reside in a 1-galaxy $\Gamma^{1}$ different 
from $\Gamma_{0}^{1}$.  $\Box$

\section{When $^{*}\!G^{1}$ Has a 1-Hypernode Not in Its Principal Galaxy}

We are at last ready to extend the results of Section 4 to the
rank 1 of transfiniteness.  The ideas are the much 
same as those of Section 4 except for the fact that the proof 
of Theorem 4.2 cannot be extended to the present case.
This is because transfinite ordinals cannot be identified 
as being even or odd.  We need another proof.

In this section $G^{1}$ is 1-wconnected and has an infinity 
of boundary 1-nodes, but $G^{1}$ need not be locally finite.
Let $\Gamma^{1}_{a}$ and $\Gamma^{1}_{b}$ be two 1-galaxies of $^{*}\!G^{1}$ 
that are different from the principal 1-galaxy $\Gamma^{1}_{0}$.
We say that $\Gamma^{1}_{a}$ {\em is closer to $\Gamma^{1}_{0}$
than is $\Gamma^{1}_{b}$} and that $\Gamma^{1}_{b}$ {\em is further away 
from $\Gamma^{1}_{0}$ than is $\Gamma^{1}_{a}$} if there are a 
${\bf y}=[y_{n}]$ in $\Gamma^{1}_{a}$ and a ${\bf z}=[z_{n}]$ in 
$\Gamma^{1}_{b}$ such that, for some ${\bf x}=[x_{n}]$ in 
$\Gamma^{1}_{0}$ and for 
every $m\in\N$,
\[ \{n\!: d(z_{n},x_{n})-d(y_{n},x_{n})\;\geq\omega\cdot m\}\;\in\;{\cal F}. \]
(The ranks of ${\bf x}$, ${\bf y}$, and ${\bf z}$ may now be either 0 or 1.)

Any set of 1-galaxies for which every two of them, say, 
$\Gamma^{1}_{a}$ and $\Gamma^{1}_{b}$ satisfy these 
conditions will be said to be 
{\em totally ordered according to their closeness
to} $\Gamma^{1}_{0}$.  Here, too, the conditions for a total 
ordering are readily shown.

{\bf Lemma 11.1.}  {\em These definitions are independent of the 
representative sequences $\langle x_{n}\rangle$, $\langle y_{n}\rangle$,
and $\langle z_{n}\rangle$ chosen for ${\bf x}$, ${\bf y}$, and 
${\bf z}$.}

The proof of this lemma is the same as that of Lemma 4.1
except that 
the rank 0 is replaced by the transfinite rank 1.  For instance, 
the natural number $m$ is now replaced by the ordinal
$\omega\cdot m$. 

{\bf Theorem 11.2.}  {\em If $^{*}\!G^{1}$ has a hypernode ${\bf v}=[v_{n}]$ 
(of either rank 0 or rank 1) that is not in its 
principal 1-galaxy $\Gamma^{1}_{0}$, then there exists a two-way 
infinite sequence of 1-galaxies totally ordered according to
their closeness to $\Gamma^{1}_{0}$ and with $\bf v$ 
being in one of those galaxies.}

{\bf Proof.} In this proof, we use the fact that between any two nodes
in a 1-graph there exists a geodesic walk terminating at those nodes;
that is, the length of the walk is equal to the wdistance 
between those nodes.  This is a consequence of the facts 
that the walks terminating at those nodes have ordinal lengths
and that any set of ordinals is well-ordered and thus has a least ordinal.
That least ordinal must be the length of at least one walk terminating 
at those nodes, for otherwise the minimum of the 
walk-lengths would be larger.

As before, let ${\bf x}=[\langle x,x,x,\ldots\rangle]$
be a standard hypernode in $\Gamma_{0}^{1}$.
${\bf x}$ can be of either rank 0 or 1.  Since $\bf v$ is not 
in $\Gamma_{0}^{1}$, we have that for each $m\in\N$
\begin{equation}
\{n\!: d(x,v_{n})\,>\, \omega\cdot m\}\;\in\;{\cal F}  \label{11.1}
\end{equation}
For each $n\in\N$, if $d(x,v_{n})< \omega\cdot 6$, set $u_{n}=x$, but,
if $d(x,v_{n})\geq \omega\cdot 6$, choose $u_{n}$ such that 
\begin{equation}
d(x,v_{n})\;\leq\; d(x,u_{n})\cdot 3\;\leq\; d(x,v_{n})\cdot 2  \label{11.2}
\end{equation}
That the latter can be done can be seen as follows.

Choose a geodesic 1-walk $W^{1}$ terminating at $x$ and $v_{n}$.
Remember that $W^{1}$, as given by (\ref{6.2}),
is incident to each of its nonterminal
1-nodes through at least one 0-tip, as was 
asserted by Condition 3 of the definition of $W^{1}$.
Moreover, the transition through each 0-tip contributes $\omega$
to the length of $W^{1}$.  Upon tracing $W^{1}$ from $x$
toward $v_{n}$, we must encounter at least two 1-nodes, 
both of which are neither closer to $x$ by one-third of the number
of 0-tips traversed by $W^{1}$ nor 
further away from $x$ by two-thirds of the number of 0-tips 
traversed by $W^{1}$. A node on $W^{1}$ between those two 1-nodes can be
chosen as $u_{n}$.

Suppose there is a $k\in\N$ such that $\{n\!: d(x,u_{n})\leq \omega\cdot k\}
\in {\cal F}$. By the left-hand inequality of (\ref{11.2}), 
\[ \{n\!: d(x,v_{n})\,\leq\, (\omega\cdot k)\cdot 3\}\;\supseteq\;
\{n\!: d(x,u_{n})\,\leq\, \omega\cdot k\}\;\in\;{\cal F}. \]
Hence, the left-hand set is a member of $\cal F$, in 
contradiction to (\ref{11.1}).  (These sets cannot both be 
in the ultrafilter $\cal F$.)  Therefore, ${\bf u}=[u_{n}]$
satisfies (\ref{11.1}) for every $m\in\N$ when $v_{n}$ 
is replaced by $u_{n}$;  that is, $\bf u$ is in a galaxy 
different
from the principal 1-galaxy $\Gamma_{0}^{1}$.

Furthermore, by the right-hand inequality of (\ref{11.2}),
\[ d(x,u_{n})\;\leq\; (d(x,v_{n})\,-\,d(x,u_{n}))\cdot 2. \]
Suppose there exists a $j\in\N$ such that 
\[ \{n\!: d(x,v_{n})\,-\, d(x,u_{n})\,\leq\, \omega\cdot j\}\;\in\;{\cal F}. \]
Then,
\[ \{n\!: d(x,u_{n})\,\leq\,(\omega\cdot j)\cdot 2\}\;\supseteq\;
\{n\!: d(x,v_{n})\,-\, d(x,u_{n})\,\leq\, \omega\cdot j\}\;\in\;{\cal F} \]
So, the left-hand set is in $\cal F$, in contradiction to our previous 
conclusion that $\bf u$ satisfies (\ref{11.1}) 
with $v_{n}$ replaced by $u_{n}$.  We can conclude that 
$\bf u$ and $\bf v$ are in different 1-galaxies $\Gamma_{a}^{1}$
and $\Gamma_{b}^{1}$ respectively, with $\Gamma_{a}^{1}$ closer
to $\Gamma_{0}^{1}$ than is $\Gamma_{b}^{1}$.

We can now repeat this argument with $\Gamma_{b}^{1}$ replaced
by $\Gamma_{a}^{1}$ and with ${\bf u}=[u_{n}]$ playing the role that
${\bf v}=[v_{n}]$ played
to find still another galaxy $\Gamma_{a}^{1}\prime$ different from 
$\Gamma_{0}^{1}$ and closer to $\Gamma_{0}^{1}$ than is $\Gamma_{a}^{1}$.  
Continual repetitions yield an infinite sequence of galaxies 
indexed by, say, the negative integers and totally ordered
by their closeness to $\Gamma_{0}^{1}$.

The rest of the proof continues just like the argument for 
Theorem 4.2 that establishes a sequence of 1-galaxies progressively
further away from $\Gamma_{0}^{1}$ than is $\Gamma_{b}^{1}$.
In this case, the natural number $m$ is replaced by the ordinal
$\omega\cdot m$;  also, the last 
$n$ in (\ref{4.2}) is replaced by $\omega\cdot n$.
$\Box$

Finally, by mimicking the proof of Theorem 4.3, we can prove

{\bf Theorem 11.3.} {\em Under the hypothesis of Theorem 11.2, the set of 
1-galaxies of $^{*}\!G^{1}$ is partially ordered according
to the closeness of the 1-galaxies to $\Gamma^{1}_{0}$.}

\section{Extensions to Higher Ranks of Transfiniteness}

The extension of these results to the enlargements of 
transfinite graphs of any natural-number rank is quite similar
to what we have presented.  The ideas are the same, but the 
notations and the details of the arguments are
somewhat more complicated.
Moreover, further complications arise with the extension to the 
arrow rank $\vec{\omega}$ of transfiniteness.
Extensions to still higher ranks then proceed in much the same way.
All this is explicated in the technical report \cite{gal2}, 
which can also be found in the archival web site, www.arxiv.org, under
Mathematics, Zemanian.

\end{document}